\documentclass[11pt]{amsart}
\usepackage{amscd, amsmath, amsthm, amssymb, xy, xypic}
\newtheorem{theorem}{Theorem}[section]

\newtheorem{lem}[theorem]{Lemma}
\newtheorem{prop}[theorem]{Proposition}

\theoremstyle{defi}     % italic or bold etc.
\newtheorem{defi}[theorem]{Definition}

\theoremstyle{remark}
\newtheorem{remark}[theorem]{Remark}
\numberwithin{equation}{section}

\def \hd #1 {\bfseries #1  \mdseries}
\def \italic #1 {\bfseries \it #1 \rm \mdseries}
\def \ra {\rightarrow}
\def \cen #1 { \begin{center} #1 \end{center}}

\def \mbc {\mathbb C}
\def \mbp {\mathbb P}
\def \mbv {\mathbb V}
\def \mbq  {\mathbb {Q}}

\def \mco  {\mathcal {O}}
\def \mci  {\mathcal {I}}
\def \mca  {\mathcal {A}}

\def \mcf  {\mathcal {F}}
\def \Q {${\mathbb {Q}}\,$}

\def \Pic {{\rm{Pic}}}
\def \rk {{\rm{rk}}}
\def \im {{\rm{im}}}

\def \Sing {{\rm{Sing}}}

\def \Spec {{\rm{Spec}}}
\def \id {{\rm{id}}}
\def \Cl {{\rm{Cl}}}
\def \pr {{\rm{pr}}}

\begin{document}
\title[Calabi--Yau coverings over  singular varieties]
{Calabi--Yau coverings over some singular varieties and new Calabi-Yau 3-folds with Picard number one}

\author[N.-H. Lee]{Nam-Hoon Lee }
\address{School of Mathematics, Korea Institute for Advanced Study, Dongdaemun-gu, Seoul 130-722, Korea }
\email{nhlee@kias.re.kr}
%\thanks{Partially supported by NSF grant DMS-0093542}
\subjclass[2000]{14J32, 14J45, 14D06}
\begin{abstract}
This note is a report on the observation that some singular varieties admit Calabi--Yau coverings.
As an application,  we construct 18 new  Calabi--Yau 3-folds with Picard number one that have some interesting properties.

\end{abstract}
\maketitle
% \tableofcontents
\setcounter{section}{-1}
\section{Introduction}
A \italic{Calabi--Yau manifold} is a compact K\"ahler manifold with trivial canonical class such that the intermediate cohomologies of its structure sheaf are all trivial  ($h^i (X, \mco_X) =0$ for $0 < i < \dim (X)$).
One handy way of constructing Calabi--Yau manifolds is by taking coverings of some smooth varieties such that some multiples of their anticanonical class have global sections. Indeed many of known examples of Calabi--Yau 3-folds with Picard number one are constructed in this way (see, for example, Table 1 in \cite{EnSt}). In this note we show that singular varieties with some cyclic singularities also admit Calabi--Yau manifolds as their coverings (Theorem \ref{2.1}). We give some formula for calculating their invariants by using degeneration method (Theorem \ref{inv}, Theorem \ref{tinv}).

In his beautiful papers (\cite{Ta1}, \cite{Ta2}), H. Takagi classified possible invariants of certain \Q-Fano 3-folds of Gorenstein index 2 and constructed some exotic examples of \Q-Fano 3-folds.
We apply our theorem to construct Calabi--Yau 3-folds which
are double coverings of  Takagi's
\Q-Fano 3-folds.
It turns out that 18 of them  are   new Calabi--Yau 3-folds with \emph{Picard number one} (Table \ref{cyt}). Although a huge number of Calabi--Yau 3-folds have been constructed, those with Picard number one are still quite rare (for example, see Table 1 in \cite{EnSt}). Note that they are primitive and play an important role in the moduli spaces of all Calabi--Yau 3-folds (\cite{Gr}).
 We  show that some of them are connected by projective flat deformation although they are of different topological types (Theorem \ref{hil}).
It is interesting that three of them have the invariants which were
predicted by C. van Enckevort and D. van Straten in
their paper (\cite{EnSt}).

Let us recall a notation for certain singularities.
Let  $a_1, \cdots, a_n$ be integers and
let $x_1, \cdots, x_n$ be coordinates on $\mbc^n$. Suppose that
the cyclic group $G$ acts on $\mbc^n$ via
\cen{ $x_i \mapsto \varepsilon^{a_i}x_i $, for all $i$,}
where $\varepsilon$ is a primitive $r^\text{th}$ root of unity for some positive integer $r$. A singularity
$q \in X$ is called a {{quotient singularity of type}}
$\frac{1}{r} (a_1, \cdots, a_n)$ if there is a neighborhood of $q$ that
is isomorphic to a neighborhood  of $(0, \cdots, 0)$ in $\mbc^n / G$.

Consider a simple example of covering.  Let
$$X = \{x^{10} + y^{10} + z^{10} + w^5 = t^2 \} \subset \mbp (1, 1, 1, 2, 5),$$
where $x$, $y$, $z$, $w$ and $t$ are homogeneous coordinates of weights 1,
1, 1, 2 and 5 respectively. Then $X$ is a Calabi--Yau 3-fold. Define a projection
$$\pi: X \ra  \mbp  (1, 1, 1, 2)$$
by
$$ (x, y, z, w, t) \mapsto  (x, y, z, w).$$
Note that $\mbp (1, 1, 1, 2)$ has a
singularity of type $\frac{1}{2} (1, 1, 1)$ at $p =  (0, 0, 0, 1)$.
Let $S = \{x^{10} + y^{10} + z^{10} + w^5 = 0 \} \subset \mbp (1, 1, 1, 2)$, then $S$ is smooth. It is
easy to see
that $X$ is a double covering of $\mbp (1, 1, 1, 2)$ with the branch locus
$$S \cup \{ p \}.$$
Over $p$, the map $\pi$ locally looks as the quotient map
$$\mbc^3 \ra \mbc^3 / \sim,$$
where $ (x, y, z) \sim  (-x, -y, -z)$.

In the following section, we prove that a variety with  singularities of type $\frac{1}{2} (1, \cdots, 1)$
admits a Calabi--Yau double
covering when the dimension is odd.

\medskip
{\it Acknowledgments}

\medskip
 Prof. Hiromichi Takagi kindly calculated the topological Euler numbers of his examples of \Q-Fano 3-folds. The author got some advice on construction of `quasi-line bundle' in  Section \ref{hodge} from Prof. J\'anos Koll\'ar.
The author is also thankful to the referee for pointing out several typos and unclear arguments in the previous version. Finally the author would like to express his sincere thanks to  Prof. Igor V. Dolgachev for many suggestions, discussions and encouragement.

\section{Existence of Calabi--Yau covering}

For  a Cartier divisor $D$  on a variety $X$,
let $\mbv (D)$ denote  the total
space of the line bundle, corresponding to $D$ and  $\pr: \mbv (D) \ra X$ be the projection.
 Note that
an element of $H^0 (X, \mco (D))$ can be regarded
as a map from  $X$ to  $\mbv (D)$.

For any integer $r$, we can define a natural map, $\mbv (D) \ra
\mbv (rD)$,  by
$$ (x, t) \mapsto  (x, t^r)$$
in local coordinates,  where $t \in \pr^{-1} (x)$.
 We denote singularities of types
$\frac{1}{2} (1, 1)$, $\frac{1}{2} (1, 1, 1)$, $\cdots$ by
$\frac{1}{2} (1^{[2]})$, $\frac{1}{2} (1^{[3]})$, $\cdots$ respectively.
Note that $\frac{1}{2}(1,1)$ is the ordinary double point singularity on a surface.

Let an $n$-dimensional variety $Y$ have a singularity of type
$\frac{1}{2} (1^{[n]})$ at $p \in Y$.
Let $f: \widetilde Y \ra Y$ be the blow-up at $p$.  Then the singularity is
resolved and  the exceptional divisor $E$ is a
copy of $\mathbb P^{n-1}$ with normal bundle $\mco_E(-2)$.
Note that
$$2 K_{\widetilde Y} = f^*  (2 K_Y) +  (n-2) E.$$

\begin{theorem}\label{2.1} Let $Y$ be a projective  $n$-dimensional
variety with singularities of type only $\frac{1}{2}  (1^{[n]})$. Suppose that $n$ is odd and that
the complete linear system $|-2K_Y|$ contains a smooth $(n-1)$-fold $S$. Then
there is a smooth projective $n$-fold $X$ with $K_X=0$ that is a double covering of $Y$ with
the
branch locus
\cen{$ S \cup {\rm{Sing}} (Y)$.}
Furthermore assume that
\cen{$h^i  (Y,  \mco_Y) = 0$ for $ 0< i < n$,}
then $X$ is a Calabi--Yau $n$-fold.
\end{theorem}
\begin{proof} For simplicity, let ${\rm{Sing}} (Y) = \{ p \}$ be composed of a single point. The proof for general case is similar.
 Let $f:
\widetilde Y \ra Y$
be the blow-up of $Y$ at $p$. Then $\widetilde Y$ is smooth. As
mentioned before, $2 K_{\widetilde Y} = f^* ( 2 K_Y ) +  (n-2) E$, where
$E$ is
the
exceptional divisor. We can set $n = 2 m + 1$ for some positive integer $m$ because $n$ is odd. Then we have
\cen{$f^* (-2K_Y) + E =  - 2 K_{\widetilde Y} +  (n-1)E =
2 (-K_{\widetilde Y} + m E)$.}
Let $ D = -K_{\widetilde Y} + m E $, then $2D = f^* (-2K_Y) + E$. Note
that $S' =
f^* (S)$ is smooth  and disjoint from $E$. Clearly
there exists a
global section $\rho \in H^0 (\widetilde Y, \mco_Y (2D))$ such that
${\rm{div}} (\rho) = S' + E$.  Note that $\rho$ is also interpreted as a
map from $Y$ to $\mbv (2D)$. Let
$\widetilde X = \phi^{-1} ({\rm{im}} \,\rho)$, where $\phi:\mbv (D)  \ra
\mbv (2D)$ is the natural map. Then $\widetilde X$ is smooth.
Let $\tilde \pi: \widetilde X \ra \widetilde Y$ be the restriction of the projection $\mbv(D)\ra \widetilde Y$ to $\widetilde X$, then $\tilde \pi$ is
a double covering map.
The canonical class of $\widetilde X$ is:
\begin{equation}\label{30}
K_{\widetilde X} =\tilde \pi^* (K_{\widetilde Y} + D) = \tilde
\pi^* (mE) = 2mF,
\end{equation}
 where $F = \tilde \pi^{-1} (E)$. Note that $F$ is isomorphic to $\mathbb
P^{n-1}$
because
$E$ is.
Let $H$ be the unique ample generator of ${\rm{Pic}} (F)$, then
$$-n H = K_F =  (K_{\widetilde X} + F) | _F                = (2mF + F)|_F
          = (2m + 1 ) F |_F                =n F|_F.$$
Since ${\rm{Pic}} (F)$ is torsion-free, $F|_F = -H $. Therefore $F$ can
be smoothly
contracted to a point. Let $g: \widetilde X \ra X$ be the contraction.
Note that $X$
is
smooth. Since
\cen{$K_{\widetilde X} = g^* (K_X) +  (n-1)F$,}
we have
\cen{ $g^* (K_X) = K_{\widetilde X} -  (n-1)F = 2mF-  (n-1)F = 0 $.}
Because the map $g^*:\rm{Pic} (X) \ra \rm{Pic} (\widetilde X)$ is injective, the canonical class of $X$ is trivial: $K_X =
0$.
Furthermore since $Y$ is projective, so are $\widetilde Y$, $\widetilde X$ and $X$.
So the first assertion is proved.
Now we assume that
\cen{$h^i  (Y,  \mco_Y) = 0$ for $ 0< i < n$}
and show that
{$h^i(X, \mco_X) =0$ for $0 < i < n$.}
Firstly we note that

\cen{$h^i  (\widetilde X, \mco_{\widetilde X}) = 0$ for $ 0< i < n$.}

Consider the following exact sequence,
$$ 0 \ra \mco_{\widetilde Y} (kE) \ra \mco_{\widetilde Y} ( (k+1)E)  \ra
\mco_E ( (k+1)E) \ra 0,$$
where $k$ is a non-negative integer. Note that
$$\mco_E ( (k+1)E)=\mco_E (-2 (k+1)).$$
We have an exact sequence
\cen{$H^{i-1} (E, \mco_E (-2 (k+1))) \ra H^i (\widetilde Y,
\mco_{\widetilde Y} (kE)) \ra H^i (\widetilde Y, \mco_{\widetilde
Y} ( (k+1)E)) \ra H^{i} (E, \mco_E (-2 (k+1)))$.}
Note that
$$H^{i-1} (E, \mco_E (-2 (k+1))) =0$$
and
$$H^{i} (E, \mco_E (-2 (k+1))) =0$$
for $0 < i <n $, $k=0,1, \cdots, m-1$, which implies that
$$ h^i (\widetilde Y, \mco_{\widetilde Y} ( (k+1)E)) =
h^i (\widetilde Y, \mco_{\widetilde Y} (kE)).$$
 Therefore we have
$$h^i (\widetilde Y, \mco_{\widetilde Y} (mE)) =h^i (\widetilde Y,
\mco_{\widetilde Y} ( (m-1)E)) = \cdots =h^i (\widetilde Y,
\mco_{\widetilde Y}) = 0 $$
for $0<i<n$.
Note that the map $\tilde \pi$ is finite and $$\tilde \pi_*
 (\mco_{\tilde{X}}) \simeq \mco_{\widetilde Y} \oplus \mco_{\widetilde
Y}  (-D) = \mco_{\widetilde Y} \oplus \mco_{\widetilde Y} (K_{\widetilde
Y} - mE).$$
Finally we have
\begin{align*}
H^i (\widetilde X, \mco_{\widetilde X}) &\simeq H^i (\widetilde Y, \tilde
\pi_* (\mco_{\tilde{X}})) \\
&\simeq  H^i (\widetilde Y, \mco_{\widetilde Y}) \oplus H^i (\widetilde
Y,  \mco_{\widetilde Y} (K_{\widetilde Y} - mE))\\
&\simeq H^i (\widetilde Y, \mco_{\widetilde Y}) \oplus H^i (\widetilde Y,
  \mco_{\widetilde Y} ( mE))^*\\
&= 0
\end{align*}
for $0< i < n$, which leads to:
$$h^i (X, \mco_X) =h^i(\widetilde X, \mco_{\widetilde X}) = 0$$
 for $0<i<n$.
In conclusion, $X$
is a Calabi--Yau $n$-fold. There is a rational map $\pi:X \dasharrow Y$ such
that the following diagram (Figure \ref{dia}) is commutative.
\begin{figure}[h]
$$
\xymatrix{
\widetilde X \ar[r]^{\tilde \pi} \ar[d]_g    &   \widetilde Y \ar[d]^f\\
X     \ar@{-->}[r]^{\pi}   &Y
}
$$\caption{}\label{dia}
\end{figure}
But it is not hard  to see that $\pi$ is actually a double covering
morphism with
branch locus $S \cup \{p\}$.
\end{proof}

\section{Hodge numbers of the Calabi--Yau  covering}\label{hodge}

For a double covering with dimension higher than two, it is  a non-trivial task to calculate the Hodge numbers even in the case that the base of the covering is smooth. They are calculated for some special cases (see, for example, \cite{Cl2}, \cite{Cy} and \cite{CySt}).
In this section, we give a formula for $h^{1,1}(X)$ when the complete linear system $|-K_X|$ contains a variety with a certain mild singularities.
In three-dimensional case, the other Hodge numbers of a Calabi--Yau 3-fold are determined by its topological Euler number, which is easily calculable.

\begin{theorem}\label{inv} Let $X$ be the Calabi--Yau
double covering of $Y$ in Theorem \ref{2.1}.
Suppose that  the linear system $|-K_Y|$ contains a
variety $D$ such
that:
\begin{enumerate}
\item The variety $D$ has its singularities at those of $Y$, \textit{i.e.} $\Sing(D)=\Sing(Y)$.
\item The singularity types of  ${\rm{Sing}} (D)$ are all $\frac{1}{2} (1^{[n-1]})$ and $C=S \cap D$ is smooth.
\item $h^1(\mco_D)=0$.
\end{enumerate}
Then $h^{1,1}(X)=2h^{2}(Y)-k$, where $k$ is the dimension of the image of the restriction map $H^2(Y, \mbq) \ra H^2(D, \mbq)$.
\end{theorem}

To prove this theorem, we make a semistable degeneration of $X$ and use the well-known Clemens--Schmid exact sequence.

Note that the anticanonical divisor $-K_Y$ is not Cartier unless $\Sing (Y) = \emptyset$. So there is no line bundle corresponding to $-K_Y$ a priori. Let $\mcf=\mco_Y(K_Y)$ and
$$\mcf^{[i]}=(\mcf^{\otimes i})^{**}$$
 be the double dual of $\mcf^{\otimes i}$.
Consider
a sheaf of algebras on $Y$,
$$\mca = \bigoplus_{i=0}^{\infty} \mcf^{[i]},$$
where we use the multiplication
$$\mcf^{[i]} \otimes \mcf^{[j]} \simeq \mcf^{[i+j]}.$$
Let $\mbv (-K_Y) = \Spec_Y (\mca)$.
Note that it is a generalization of the construction of a line bundle for a
Cartier divisor  (\textit{cf.} Definition 5 in \cite{Ko}).  Let
$$\pr: \mbv (-K_Y)  \ra Y$$
be the projection. Let
$$\mci = \bigoplus_{i=1}^{\infty}\mcf^{[i]}.$$
The map
$$\mca \ra \mca/\mci = \mco_Y$$ induces the embedding
\begin{align*}
\Spec_Y(\mco_Y) \hookrightarrow \Spec_Y(\mca),
\end{align*}
\textit{i.e.}
\begin{align}
i: Y \hookrightarrow  \mbv (-K_Y),\label{sub}
\end{align}
under which we
consider $Y$ as a subvariety of $\mbv (-K_Y)$.
Note that $\mbv(-2K_Y)=\Spec(\mca^{(2)})$, where
$$\mca^{(2)} =  \bigoplus_{i=0}^{\infty} \mcf^{[2i]}.$$
The map
$$\Spec_Y(\mca) \ra \Spec_Y(\mca^{(2)})$$
is denoted by
\begin{align}
\psi: \mbv (-K_Y) \ra \mbv (-2K_Y).\label{map}
\end{align}

Let $Y^*=Y \setminus \Sing(Y)$, then we have
$$\pr^{-1} (Y^*) \simeq \mbv (-K_{Y^*}).$$
Note that $-K_{Y^*}$ is Cartier on $Y^*$.

This kind of construction was studied in a more general setting under the name of `Seifert $G_m$-bundle' (\cite{Ko}).
See Section 2 of \cite{Ko} for relevant discussion.
The variety $\mbv(-K_Y)$ has its singularities at those of $Y$
and  the singularities in $\mbv (-K_Y)$ are all of type $\frac{1}{2} ( 1^{[n+1]})$, where we consider $Y$ as a subvariety of $\mbv(-K_Y)$ (equation \ref{sub}). Over
singular points of $Y$, the projection  $\pr: \mbv (-K_Y)  \ra Y$ locally looks
like the following map:
$$\phi: V_n \ra U_n,$$
where \begin{align*}
V_n &= \mbc^{n+1} /  (x_1, \cdots, x_{n}, t) \sim  (-x_1, \cdots, -x_{n}, -t) \\
U_n &=  \mbc^{n} /  (x_1, \cdots, x_{n}) \sim  (-x_1, \cdots, -x_{n})
\end{align*}
and the projection $\phi$ maps:
$$ (x_1, \cdots, x_{n}, t) \mapsto  (x_1, \cdots,
x_{n}).$$

We call $\mbv (-K_Y)$ with the projection $\pr: \mbv (-K_Y) \ra Y$ a `quasi-line
bundle' corresponding to $-K_Y$.

A section $\rho$ of this quasi-line bundle   is a map $\rho: Y \ra \mbv (-K_Y)$ such that
$\pr \circ
\rho = \id_Y$.
Note that the map $ \pr|_{{\rm{im}}\,\rho}: {\rm{im}} \,\rho \ra Y$ is an isomorphism.
We obtain the double coverings in Theorem \ref{2.1}  in this setup.

\begin{prop}   Let $Y$ be a projective  $n$-dimensional
variety $Y$ with singularities of type only $\frac{1}{2}  (1^{[n]})$. Suppose that $n$ is odd and that
the complete linear system $|-2K_Y|$ contains a smooth $ (n-1)$-fold $S$. Then
there is a smooth projective $n$-fold $X$ with $K_X=0$ that is a double covering of $Y$ with
the
branch locus
\cen{$ S \cup {\rm{Sing}} (Y)$.}
\end{prop}
\begin{proof}
Let $\rho$ be a section of
$\mco_Y (-2K_Y)$ such that ${\rm{div}} (\rho) = S$. Let $X =
\psi^{-1} ({\rm{im}} \,\rho)$(see equation \ref{map} for $\psi$).
We want to show that $X$ is smooth.
Firstly, the open set
$ \psi^{-1} (
\rho (Y^*))$ of $X$ is smooth. Let us look at what happens near  $ \psi^{-1} (
\rho ({\rm{Sing}} (Y)))$. Locally it looks like:
$$\{ (x_1, \cdots, x_n, t) \in V_n \big| t^2 =  g (x_1, \cdots, x_n) \},$$
where $g:U_n \ra \mbc$ is a function such that $g (0, \cdots, 0) \neq
0$ ($\because$ $S$ is smooth, so disjoint from ${\rm{Sing}} (Y)$). By
direct calculation, one
can show that this variety is smooth. Therefore $X$ is smooth. The map $\pi =  \pr|_X: X \ra Y$ is a double
covering, branched along $S \cup {\rm{Sing}} (Y)$. Note that
$$K_{\pi^{-1} (Y^*)} = \pi|_{\pi^{-1} (Y^*)}^* (K_{Y^*} +  (-K_{Y^*})) = 0$$
Since the set
$$X \setminus \pi^{-1} (Y^*) = \pi^{-1} (\Sing (Y))$$
 is  finite, we have $K_X = 0$.
\end{proof}
It is not hard to see that the Calabi--Yau manifold in the above theorem is the same Calabi--Yau manifold in Theorem \ref{2.1} for given $Y$.

Now we construct a degeneration of $X$, assuming that the linear system
$|-K_Y|$ contains a variety $D$ such that ${\rm{Sing}} (D) =
{\rm{Sing}} (Y)$, the
singularity types of ${\rm{Sing}} (D)$ are all $\frac{1}{2} (1^{[n-1]})$
and $D \cap S$
is  smooth (the condition in Theorem \ref{inv}).

Let $\alpha$ be a section of
$\mco_{Y^*} (- K_{Y^*})$
such that
$$D\cap Y^* = {\rm{div}} (\alpha).$$
It is not
hard to see that there
is a unique extension of $\alpha$ to a section $\delta:Y \ra \mbv (-K_Y)$
of the quasi-line bundle. Let
$-\delta$ be the extension of $-\alpha$.  We also set $Y_1 = {\rm{im}} (\delta)$
and $Y_2 =
{\rm{im}} (-\delta)$. Then $Y_i$ is a copy of $Y$.
If we regard $Y$ as a subvariety of $\mbv (-K_Y)$ (equation \ref{sub}), we have
$$Y_1 \cap Y_2 = Y_i \cap Y = D$$
and
$$X \cap Y = S.$$
We have a section $\alpha^2$ of $\mco_{Y^*} (-2K_{Y^*})$ and it has a
unique extension to a
section of $\mco_Y(-2K_Y)$. Let us denote it by $\delta^2$.
Note that $\delta ^2 + t \rho$ is a section of $\mco_Y(-2K_Y)$ for $t \in \mbc$.
Let $\Delta$ be a small open disk which is centered at the origin
in $\mbc$. We define a variety:
$$Z = \bigcup_{t \in \Delta} \psi^{-1} ({\rm{im}} ( \delta ^2 + t \rho)
) \subset \mbv(-K_Y)
\times \Delta.$$
Let
$$f: Z \ra \Delta$$
be the projection
and $Z_t = f^{-1} (t)$. Then the central fiber is:
$$Z_0 = Y_1 \cup Y_2.$$
and a generic fiber ($t \neq 0$)
$Z_t$ is
a Calabi--Yau $n$-fold which is a deformation of $X$   (Note that $X =
\psi^{-1} ({\rm{im}} ( \rho) ) = f^{-1} (\infty)$).

By direct calculation, we have:

\begin{prop} If $\Delta$ is sufficiently small, then
$${\rm{Sing}} (Z) = \left ( (D \cap S) \sqcup {\rm{Sing}} (Y) \right) \times
\{0\} .$$
\end{prop}
Along  the smooth $C = D \cap S$, $Z$ is locally the product of an $ (n-2)$-fold and
a three-dimensional ordinary double point singularity.  If
we  blow up $Z$ along $C$, then the exceptional locus is a $\mathbb P^1
\times \mathbb
P^1$-bundle over $C$. It is a usual procedure to contract one of the
ruling of the
bundle smoothly to get $Z'$ (see, for example, Proposition II.1. in \cite{Lee}).
We can choose the ruling of the contraction such that
$$Z'_0 = Y'_1 \cup Y'_2,$$ where $Y'_1$ is isomorphic to $Y_1$ and
$Y'_2$ is the
blow-up
of $Y_2$ along $C$. Now the remaining singularities of $Z$ are those
which correspond to ${\rm{Sing}} (Y) \times \{0\}$ and they are all singularities of type
 $\frac{1}{2} (1^{[n+1]})$.
In summary, we have:
\begin{prop}\label{p4.6} Let $X$ be the Calabi--Yau
double covering of $Y$ in Theorem \ref{2.1}. Assume the existence of $D$ as in Theorem \ref{inv}.
 Then we can construct a
degeneration $Z'
\ra \Delta$ such that:
\begin{enumerate}
\item A generic fiber is a deformation of the Calabi--Yau double
covering $X$ of $Y$.
\item The central fiber is a union of $Y'_1$ and $Y'_2$ such that $Y'_1$
is a copy of $Y$ and $Y'_2$ is the blow-up of $Y$ along $D\cap S$.
$Y'_1 \cap Y'_2 = D$.
\item The total space $Z'$ has its singularity at $\Sing(Y'_1) = \Sing(Y'_2)$ and the singularity types  are all
$\frac{1}{2} (1^{[n+1]})$.
\end{enumerate}
\end{prop}

By analyzing the well-known Clemens--Schmid exact sequence(\cite{Cl1}), we have the following lemma (Theorem III.1 in \cite{Lee}):
\begin{lem} \label{clemens}Let $W_0$ be central compact simple normal crossing fiber in a semistable degeneration and $W_t$ ($t\neq0$) be the smooth compact K\"ahler fiber with $h^2(\mco_{W_t})=0$. Then we have
$$h^{1,1}(W_t)=h^2(W_t) = h^2(W_0) - r + 1,$$
where $r$ is the number of components of $W_0$ and $h^2(W_0)=\dim_Q H^2(W_0, \mbq)$, etc.
\end{lem}
For detailed proof, we refer to Chapter III of \cite{Lee}.

Now we are ready to prove Theorem \ref{inv}.
\begin{proof}[proof of Theorem \ref{inv}]
For simplicity, let $\Sing(Y)$ be a single point. The proof for the general case will be given later. Since the Hodge numbers are invariant under deformation, $h^{1,1}(X)$ is the same with that of a generic fiber of the degeneration $Z' \ra \Delta$ in Proposition \ref{p4.6}. Since $Z'$ has a singular point, the degeneration is not semistable. Let $W \ra Z'$ be the blow-up at the singular point. Then we have another degeneration, $W \ra \Delta$: the composite of $W\ra Z'$ and $Z' \ra \Delta$.
Note that the generic fiber is not affected and the central fiber  $W_0=V_1 \cup V_2 \cup E$ is a reduced normal crossing variety, where $V_i$ is the proper transformation of $Y'_i$ and $E$ is the exceptional divisor.
So $W \ra \Delta$ is a  semistable degeneration. Note that $\widetilde D := V_1 \cap V_2$ is a smooth $(n-1)$-fold that is the proper transformation of $Y'_1 \cap Y_2' \simeq D$. Let $H_i := V_i \cap E$, then it is a copy of $\mbp^{n-1}$.  Note also that $H_i \cap \widetilde D$ and $H_1 \cap H_2$ are  copies of $\mbp^{n-2}$.
Recall that $h^{2}(\mco_{W_t})=h^2(\mco_{X})=0$ for $t \neq 0$. By Lemma \ref{clemens},
\begin{align}\label{14}
h^{2}(W_t)= h^2(W_0)-3+1 = h^2(W_0)-2.
\end{align}

Note that $h^1(\mco_D)$ implies $h^1(\mco_{\widetilde D})=0$. By Hodge decomposition theorem, we have $H^1(\widetilde D, \mbq)=0$.
The Mayer--Vietoris sequence gives:
\cen{$\dim H^2(H_1 \cup H_2, \mbq)=1$ and $\dim H^1(H_1 \cup H_2, \mbq)=0$.}
Let $W'_0 = V_1 \cup V_2$ and consider the following exact sequence:
$$0 \ra \mbq_{W_0} \ra \mbq_{W_0'} \oplus \mbq_{E} \ra \mbq_{H_1 \cup H_2} \ra 0$$
to derive
\begin{align}\label{seq}
0&=H^1(H_1 \cup H_2, \mbq) \ra H^2(W_0, \mbq) \ra H^2(W_0', \mbq) \oplus H^2(E, \mbq) \ra H^2(H_1 \cup H_2, \mbq) \ra.
\end{align}
Since $H^2(E, \mbq) \ra H^2(H_1 \cup H_2, \mbq)$ is not a zero map, neither is
$$H^2(W_0', \mbq) \oplus H^2(E, \mbq) \ra H^2(H_1 \cup H_2, \mbq).$$
Because $\dim H^2(H_1 \cup H_2, \mbq)=1$, the above map is actually surjective. So the sequence \eqref{seq} becomes a short exact sequence:
\begin{align}\label{shortseq}
0& \ra H^2(W_0, \mbq) \ra H^2(W_0', \mbq) \oplus H^2(E, \mbq) \ra H^2(H_1 \cup H_2, \mbq) \ra 0.
\end{align}
Accordingly we have
\begin{align}\label{11}
h^2(W_0) &= h^2(W_0') + h^2(E)-h^2(H_1 \cup H_2)=h^2(W_0').
\end{align}
But with $H^1(\widetilde D, \mbq)=0$, again  the Mayer--Vietoris sequence gives:
\begin{align}
\nonumber h^2(W_0') =& h^2(V_1) + h^2(V_2) - \dim (\im (H^2(V_1, \mbq) \oplus H^2(V_2, \mbq) \ra H^2(\widetilde D, \mbq))) \\
\nonumber =& h^2(Y_1')+1 + h^2(Y_2')+1-\dim (\im (H^2(V_1, \mbq) \oplus H^2(V_2, \mbq) \ra H^2(\widetilde D, \mbq)))\\
\nonumber=&h^2(Y)+1 + h^2(Y)+2-\left (\dim (\im (H^2(Y_1', \mbq) \oplus H^2(Y_2', \mbq) \ra H^2(D, \mbq)))+1 \right )\\
\nonumber=&h^2(Y)+1 + h^2(Y)+2- (k+1)\\
=&2h^2(Y)+2-k\label{12} .
\end{align}
So we have
\begin{align}\label{15}
h^{1,1}(X) &= h^2(X)=h^{2}(W_t)= h^2(W_0)-2=h^2(W'_0)-2=2h^2(Y)-k.
\end{align}

Let us consider general case that $m=\#\Sing(Y)$ is not necessarily one. Let $W \ra Z'$ be the blow-up at singular points. Let
$E_1, \cdots, E_m$ be the exceptional divisors and $H_{ij}=V_i \cap E_j$. We note that $E_i \cap E_j =\emptyset$. Then
the equation \eqref{14} becomes
$$h^{2}(W_t)= h^2(W_0)-(2+m)+1 = h^2(W_0)-m-1,$$
the short exact sequence \eqref{shortseq} becomes
$$0 \ra H^2(W_0, \mbq) \ra H^2(W_0', \mbq) \oplus \bigoplus_i H^2(E_i, \mbq) \ra \bigoplus_{j}H^2(H_{1j} \cup H_{2j}, \mbq) \ra 0,$$
we still have the equation \eqref{11}:
$$h^2(W_0) = h^2(W_0') + \sum_i h^2(E_i)-\sum_j h^2(H_{1j} \cup H_{2j})=h^2(W_0')+ m -m = h^2(W_0'),$$
the equation \eqref{12} becomes:
$$h^2(W_0') = 2h^2(Y)+2m+1-(k+m)=2h^2(Y)+m+1-k$$
and finally the equation \eqref{15} becomes:
$$h^{1,1}(X) = h^{2}(W_t)= h^2(W_0)-m-1=2h^2(Y)+m+1-k-m-1=2h^2(Y)-k.$$
\end{proof}

\section{Calabi--Yau 3-folds as coverings of Takagi's \Q-Fano 3-folds}

In \cite{Ta1} and \cite{Ta2}, H. Takagi classified possible invariants of some \Q-Fano 3-folds of Gorenstein index 2  and constructed some examples of \Q-Fano 3-folds.
They ($Y$'s) have the following properties (Corollary 3.4 in \cite{Ta2}):
\begin{enumerate}
\item The singularities are all isolated, type of  $\frac{1}{2} (1, 1, 1)$.
\item The class group is: $\Cl(Y) = \langle K_Y \rangle $.
\item The divisor $-2K_Y$ is   very ample.
\item A general element of $|-K_Y|$ is a singular $K3$ surface with
only ordinary
double points at ${\rm{Sing}} (Y)$. Say $D$ is such a divisor.
\end{enumerate}
The first and second properties imply that $Y$ is a \Q-Fano 3-folds.
The third property implies that the complete linear system $|-2K_Y|$ contains a smooth surface $S$ such that an intersection $C = S \cap D$ is
a smooth curve. By Theorem \ref {2.1}, there are Calabi--Yau 3-folds that are double coverings of the above \Q-Fano 3-folds. We note that
\cen{$h^2(Y)=1$, $\dim (\im (H^2(Y, \mbq) \ra H^2(D, \mbq)))=1$ and $h^1(\mco_D) = 0$.}
The product of a divisor with second Chern class $c_2(Y)$ of $Y$ is defined as follows (\textit{cf.} Definition 10.1 in page 441, \cite{Re}).

\begin{defi}Let $Y$ be a 3-dimensional projective variety with only isolated singularities of type $\frac{1}{2}(1,1,1)$.  Let $f:\widetilde Y \ra Y$ be the blow-up at $\Sing(Y)$. Note that $\widetilde Y$ is smooth.
For any divisor $L$ in $\Cl(X)$, the integer $L \cdot c_2(Y)$ is defined as
$$L \cdot c_2(Y) = f^*(L) \cdot c_2(\widetilde Y).$$
\end{defi}

Note that Takagi's \Q-Fano 3-folds satisfy all the conditions in Theorem \ref{2.1} and Theorem \ref{inv}.
We calculate the invariants of Calabi--Yau 3-folds that are double coverings of Takagi's \Q-Fano 3-folds:

\begin{theorem}\label{tinv} Let $X$ be the Calabi--Yau double covering of one of Takagi's \Q-Fano 3-folds, $Y$ with
$(-K_Y^3, N) \neq (4,4)$, where $N=\#\Sing(Y)$.
Then $X$ has the following invariants.
\begin{enumerate}
\item It has Picard number one: $\,\,\rk \Pic(X) = h^{1,1}(X)=1$. \label{f1}
\item Let $\,{\rm{Pic}}_f (X) = \langle H \rangle$, then
$$H^3 = -2 K_Y^3,$$
where $\Pic_f(X) = \Pic(X)/{\rm torsion}$ and $H$ is an ample divisor. \label{f2}
\item The product of the divisor $H$ with the second Chern class is:
$$H \cdot c_2 (X) = 48 - 3 N - 2K_Y^3,$$
where $N=\#\Sing(Y)$. \label{f3}
\item The topological Euler number is:
$$e (X) = 2 e (Y)  + 4 K_Y^3 + 2N-48.$$ \label{f4}
\item $h^{1,2}(X)=25-e(Y)-2 K_Y^3 - N$. \label{f5}
\end{enumerate}
\end{theorem}
\begin{proof}
We will freely use notations in the proof of Theorem \ref{2.1}.
By Theorem \ref{inv}, we have the assertion \eqref{f1}:
$$\rk\Pic(X)=h^{1,1}(X)=2h^{2}(Y)-\dim(\im( H^2(Y, \mbq) \ra H^2(D, \mbq)))=2-1=1.$$
Let $S''' = \pi^{-1}(S)$. Then $2S'''=\pi^*S$ as divisor classes. Let $S''' = l H$ for some integer $l$, then
$$2l H = 2S''' = \pi^*S = -2\pi^*K_Y.$$
So  $lH = -\pi^*K_Y$ and $H^3 = -2 K_Y^3/l^3$.
Let $\pi_*(H)=l' K_Y$ for some integer $l'$. Note that $\pi_*\pi^*K_Y = 2K_Y$. Then
$$2K_Y = \pi_*\pi^*K_Y = -l \pi_*H = - l l' K_Y.$$
So $l l' = -2$.
Note
$$0 < H^3 = -\pi^*(K_Y)^3/l^3= -2 K_Y^3/l^3.$$
So $l=1$ or $2$ ($\because K_Y^3 < 0$). We will exclude the case  $l=2$ at the end of the proof. Then we have $H = - \pi^*(K_Y)$ and $H^3 = -2 K_Y^3$, which is the assertion \eqref{f2}.
Let us consider the diagram (Figure \ref{dia}) in the proof of Theorem \ref{2.1}:
$$
\xymatrix{
\widetilde X \ar[r]^{\tilde \pi} \ar[d]_g    &   \widetilde Y \ar[d]^f\\\
X     \ar[r]^{\pi}   &Y
}
$$
Note that
\begin{equation}\label{60}
g^* H=- g^* \pi^*K_Y=- \tilde \pi^* f^* K_Y.
\end{equation}
Let $S'' = \tilde\pi^{-1}(S')$ and $E_1, \cdots, E_N$ be the exceptional divisor of $f:\widetilde Y \ra Y$, $F_i = \widetilde \pi^{-1}(E_i) \simeq E_i \simeq \mbp^3$ for $i=0, \cdots, N$. Note that $S''$ is isomorphic to $S'$ and that $2S''=\tilde \pi^*(S')$ as divisor classes. By adjunction,
\begin{equation}\label{21}
c_1(\widetilde Y)|_{S'} = S'|_{S'}+c_1(S')
\end{equation}
and
\begin{equation}\label{20}
c_2(\widetilde Y)\cdot {S'} = S'|_{S'} \cdot c_1(S')+c_2(S').
\end{equation}
So
\begin{align*}
e(S')=&c_2(S')\\
=&c_2(\widetilde Y)\cdot {S'} - S'|_{S'} \cdot c_1(S') \hspace{120pt}(\because \text{ equation \eqref{20}})\\
=&c_2(\widetilde Y)\cdot f^*(-2K_Y) - S'|_{S'} \cdot(c_1(\widetilde Y)|_{S'} - S'|_{S'})\,\,\,\,\,\,\,\,\,\,\,\,\,\,(\because \text{ equation }\eqref{21})\\
=&-2 c_2( Y)\cdot K_Y -S'^2\cdot(c_1(\widetilde Y) -S')        \\
=&-2 c_2( Y)\cdot K_Y-f^*(-2K_Y)^2 \cdot(f^*(-K_Y)+\frac{1}{2}(E_1+ \cdots + E_N)-f^*(-2K_Y))\\
=&-2 c_2( Y)\cdot K_Y-4f^*(K_Y)^3\\
=&48-3N-4 K_Y^3 \hspace{150pt}(\because \text{ Theorem 1.2 in \cite{Ta1}})
\end{align*}
The equation \eqref{30} implies that $c_1(\widetilde X)|_{S''}=0$ ($\because S'' \cap F_i =\emptyset$). By adjunction,
$$0=c_1(\widetilde X)|_{S''} = S''|_{S''}+c_1(S'').$$
So we have
\begin{equation}\label{61}
 S''|_{S''}  \cdot c_1(S'') =-c_1(S'')^2 = -c_1(S')^2=2f^*(K_Y)^3=2K_Y^3.
\end{equation}
We show the assertion \ref{f3}:
\begin{align*}
l H \cdot c_2(X) =&g^* (l H) \cdot c_2(\widetilde X)\\
=&-\tilde \pi^* f^* K_Y \cdot c_2(\widetilde X)\hspace{100pt}(\because \text{ equation \eqref{60}})\\
=& S'' \cdot c_2(\widetilde X)  \hspace{100pt}(\because S''= (1/2)\tilde \pi^* S' = (1/2)\tilde \pi^* f^*(-2K_Y))\\
=& S''|_{S''} \cdot c_1(S'')+c_2(S'') \hspace{85pt}(\text{By adjunction})\\
=& 2K_Y^3+e(S'')\hspace{125pt}(\because \text{ equation \eqref{61}})\\
=&2K_Y^3+48-3N-4 K_Y^3\\
=& 48 - 3 N - 2K_Y^3.
\end{align*}
We will show $l=1$ later.

We have the assertion \ref{f4}:
\begin{align*}
e(X) &= e(\widetilde X) - \sum_i(e(F_i)-1)\\
&=e(\widetilde X) - 2 N\\
&=2e(\widetilde Y) - (e(S')+\sum_i e(E_i))-2N\\
&=2e(\widetilde Y) - (48-3N-4 K_Y^3+3N)-2N\\
&=2e(\widetilde Y) - 48+4 K_Y^3 -2N\\
&=2(e(Y)+\sum_i(e(E_i)-1)) - 48+4 K_Y^3 -2N\\
&=2e(Y)+4N - 48+4 K_Y^3 -2N\\
&=2e(Y)+2N + 4K_Y^3-48.
\end{align*}
The assertion \ref{f5} follows from the fact:
$$e(X)=2(h^{1,1}(X) - h^{1,2}(X)).$$

Now we show that $l=1$. Suppose that $l=2$, then we have
\begin{align*}
&H^3 = \frac{-\pi^*(K_Y)^3}{l^3}= \frac{-K_Y^3}{4},\\
&H \cdot c_2(X) = \frac{(48 - 3 N - 2K_Y^3)}{l} =\frac{(48 - 3 N - 2K_Y^3)}{2}.
\end{align*}
By  Riemann--Roch formula, we have:
\begin{align*}
 \chi(\mco_X(H))    =& \frac{H^3}{6} + \frac{H \cdot c_2(X)}{12} \\
            =& \frac{-K_Y^3}{24} + \frac{48 - 3 N - 2K_Y^3}{24}\\
            =& 2  + \frac{-K_Y^3-N}{8}.
\end{align*}
Note that
\cen{$H^3 = \frac{-K_Y^3}{4}$ or $\chi(\mco_X(H))-2=\frac{-K_Y^3-N}{8}$}
is not an integer for
Takagi's Fano 3-folds except the case  $(-K_Y^3, N) = (4,4)$. So the number $l$ cannot be 2 and we have  $l=1$.
\end{proof}

\begin{remark}One of Takagi's \Q-Fano 3-folds has the invariant, $(-K_Y^3, N) = (4,4)$. The assertions (1) and (4) in Theorem \ref{tinv} still hold for it.
\end{remark}

Table \ref{cyt} lists the invariants of some of Takagi's \Q-Fano 3-folds and their Calabi--Yau double coverings.

\begin{table}[h]
$$\begin{array}{|c|r|c|r||c||c|r|c|r|}
\hline
\text{C.--Y. \#} &H^3 &   c_2 (X) \cdot H &{{e}} (X) &
\Leftarrow  &\text{Fano } \#   &   -K_Y^3  &   N&{{e}} (Y) \\ \hline
1&12    &   48  &-68&     &    4.2 &   6   &   4&-2    \\
2&14    &   56  &-96&    &    1.1 &   7   &   2&-12   \\
3&14    &   56  &-100&    &    4.4 &   7   &   2&-14   \\
4&15    &   54  &-76&    &    1.2, 1.3 &   7.5 &   3&-2    \\
5&15    &   54  &-84&    &    4.5 &   7.5 &   3&-6    \\
6&16    &   52  &-72&   &    4.6 &   8   &   4&0 \\
7&17    &   62  &-108&   &    1.4 &   8.5 &   1&-14   \\
8&17    &   50  &-64&    &    4.7 &   8.5 &   5&4 \\
9&18    &   60  &-84&   &    1.5 &   9   &   2&-2    \\
10&18    &   60  &-92&    &    1.6 &   9   &   2&-6    \\
11&19    &   58  &-76&    &    1.7, 1.8 &   9.5 &   3&2 \\
12&21    &   66  &-100&    &    1.9 &   10.5    &   1&-6    \\
13&21    &   66  &-104&   &    1.10    &   10.5    &   1&-8    \\
14&22    &   64  &-92&    &    1.11    &   11  &   2&-2    \\
15&25    &   70  &-100&  &    1.12    &   12.5    &   1&-2    \\
16&29    &   74  &-104&   &    1.13    &   14.5    &   1&0 \\
17&29    &   74  &-96&   &    4.8 &   14.5    &   1&4 \\
18&30    &   72  &-96&    &    1.14    &   15  &   2&4 \\ \hline
\end{array}$$
\caption{New Calabi--Yau 3-folds as double coverings of Takagi's \Q-Fano 3-folds} \label{cyt}
\end{table}

\begin{remark}[for Table \ref{cyt}]\

\begin{enumerate}
\item The rightmost four columns are for Takagi's examples. The `Fano \#' is
the index number of the \Q-Fano 3-folds in \cite{Ta1}, \cite{Ta2}.
\item The leftmost four columns are for the Calabi--Yau 3-folds which
are double coverings of  Takagi's examples.
\item H. Takagi  calculated the topological Euler numbers of his \Q-Fano 3-folds, which are not provided in \cite{Ta1} and \cite{Ta2}.
\item We only list the cases that:
    \begin{itemize}
        \item H. Takagi constructed concrete examples (see \S4. in \cite{Ta2}) and,
        \item  The invariants of the corresponding Calabi--Yau double coverings do not overlap with known examples in  Table 1 of \cite{EnSt}.
    \end{itemize}
\end{enumerate}
\end{remark}

Since the invariants:
$$ (H^3, \,H \cdot c_2 (X),\, {{e}} (X))$$
are topological, Table \ref{cyt} contains 18 new topologically different examples of Calabi--Yau 3-folds with Picard number one. Note that they can not be connected by smooth deformation.
It is interesting that some of them are, nevertheless, connected by projective flat deformation:

\begin{theorem}\label{hil} The following groups of Calabi--Yau 3-folds:
\cen{$\{ \#2, \#3 \}, \{ \#4, \#5 \}, \{ \#9, \#10 \}, \{ \#12, \#13 \}$ and  $\{ \#16, \#17 \}$}
belong to the same Hilbert schemes of some projective space,  where `$\#$' is the Calabi--Yau number in Table \ref{cyt}. Accordingly they are connected by (necessarily non-smooth) projective flat deformation.
\end{theorem}
\begin{proof} Note that the Calabi--Yau 3-folds in each group have the same values of
$$(H^3, H \cdot c_2(X)).$$
By Kodaira vanishing theorem and Riemann--Roch formula, we have:
$$\dim(H^2( \mco_X(n H))) = \chi(\mco_X(nH)) = \frac{H^3}{6}n^3 + \frac{H \cdot c_2(X)}{12} n$$
for any positive integer $n$. Note that all the Calabi--Yau 3-folds in each group have the same polynomial of $n$.
We know that $8H$ is very ample on $X$ (\cite{GaPu}). So all the Calabi--Yau 3-folds have embeddings to the same projective space $\mbp^N$ by the linear systems $|8H|$ and have the same Hilbert polynomials, where $N=\dim(H^2( \mco_X(8 H)))-1$. Therefore they belong to the same Hilbert Scheme. The final assertion follows from the fact that a Hilbert scheme is connected (\cite{Ha}).
\end{proof}

In Figure \ref{cyplot}, we plot their Hodge numbers with those of complete intersection Calabi--Yau 3-folds in weighted four-space (\cite{CaOsKa}). The topological Euler number is plotted horizontally and the number $h^{1,1} + h^{1,2}$ vertically.
 Our new Calabi--Yau 3-folds are marked by bigger dots. Note that they sit on the boundary of the plot. It would be  interesting to find the mirror partners of them.

In an interesting paper (\cite{EnSt}), C. van Enckevort and D. van Straten predicted new Calabi--Yau 3-folds
with Picard
number one.
They considered possible fourth order Picard-Fuchs equations for Calabi--Yau
3-folds with $h^{1,2} = 1$ and calculated $H^3$,
$H \cdot c_2 (X)$ and $e (X)$ for the mirror pairs whose Picard numbers are
all one.
In their list, three of the conjectural Calabi--Yau
3-folds have the values of $H^3$, $H \cdot c_2 (X)$ and $e(X)$ that show up in
the above
list. Those $ (H^3, H \cdot c_2 (X), {{e}} (X))$'s are:
\cen{$ (14, 56, -100)$, $ (21, 66, -100)$ and $ (25, -70, -100)$.}

\clearpage

\begin{figure}[h]
\input{cyplot}
\caption{}\label{cyplot}
\end{figure}

\clearpage
%%%%%%%%%%%%%% reference %%%%%%%%%%%%%%%%%

\end{document}